\documentclass[12pt,a4paper]{article}

\usepackage{color}

\usepackage{theorem,enumerate}
\usepackage{amsmath,latexsym,amssymb,amsfonts}
\usepackage{eucal}
\usepackage{mathrsfs}
\usepackage{array}
\usepackage{footmisc}
\usepackage{comment}
\usepackage{graphicx, color}

\usepackage{bm}

\newtheorem{Lem}{Lemma}[section]
\newtheorem{Thm}[Lem]{Theorem}

\newtheorem{Cor}[Lem]{Corollary}

\newtheorem{Def}[Lem]{Definition}

\title{The Alon-Tarsi number of $K_5$-minor-free graphs}
\author
{
Toshiki Abe$^{\rm a}$,
Seog-Jin Kim$^{\rm b}$,
Kenta Ozeki$^{\rm a}$
\\
{\footnotesize$^{\rm a}$Graduate school of Environmental and Information Science, Yokohama National University}\\
{\footnotesize$^{\rm b}$Department of Mathematics Education, Konkuk university}\\
}
\begin{document}
\maketitle

\begin{abstract}
In this paper, we show the following three theorems. Let $G$ be a $K_5$-minor-free graph. Then Alon-Tarsi number of $G$ is at most $5$, there exists a matching $M$ of $G$ such that the Alon-Tarsi number of $G-M$ is at most $4$, and there exists a forest $F$ such that the Alon-Tarsi number of $G-E(F)$ is at most $3$.   
\end{abstract}

\section{Introductions}
In this paper, we only deal with finite and simple graphs. A {\it $d$-defective coloring}
 of $G$ is a coloring $c:V(G)\to \mathbb{N}$ such that each color class induces a subgraph of maximum degree at most $d$.
Especially, a $0$-defective coloring is also called a {\it proper coloring} of $G$.

A {\it $k$-list assignment} of a graph $G$ is a mapping $L$ which assigns to each vertex $v$ of $G$ a set $L(v)$ of $k$ permissible colors.
Given a $k$-list assignment $L$ of $G$, a $d$-defective $L$-coloring of $G$ is a $d$-defective coloring $c$ such that $c(v)\in L(v)$ for every vertex $v$.
We say that $G$ is {\it $d$-defective $k$-choosable}  if $G$ has a $d$-defective $L$-coloring for every $k$-list assignment $L$.
Especially, we say that $G$ is {\it $k$-choosable} if $G$ is $0$-defective $k$-choosable.
The {\it choice number} $ch(G)$ is defined as the smallest integer $k$ such that $G$ is $k$-choosable.

Let $G$ be a graph and let `$<$' be an arbitrary fixed ordering of the vertices of $G$.
The \emph{graph polynomial} of $G$ is defined as 
 $$P_{G}(\bm{x})=\prod_{u\sim v,u<v}(x_u-x_v),$$
 where $u\sim v$ means that $u$ and $v$ are adjacent, and $\bm{x}=(x_v)_{v\in V(G)}$ is a vector of $|V(G)|$ variables indexed by the vertices of $G$.
It is easy to see that a mapping $c:V(G) \to \mathbb{N}$ is a proper coloring of $G$ if and only if $P_{G}(\bm{c}) \neq 0$, where $\bm{c} = \big(c(v) \big)_{v \in V(G)}$.
Therefore, to find a proper coloring of $G$ is equivalent to find an assignment of $\bm{x}$ so that $P_{G}(\bm{x}) \neq 0$.
The following theorem, which was proved by Alon and Tarsi, gives sufficient conditions for the existence of such assignments as above.

\begin{Thm}[\cite{Alon1999}]\label{cnull}(Combinatorial Nullstellensatz) Let $\mathbb{F}$ be an arbitrary field and let $f=f(x_1,x_2,\ldots,x_n)$ be a polynomial in $\mathbb{F}[x_1,x_2,\ldots,x_n]$. Suppose that the degree $\deg(f)$ of $f$ is $\sum_{i=1}^n t_i$ where each $t_i$ is a nonnegative integer, and suppose that the coefficient of $\prod_{i=1}^n x_i^{t_i}$ of $f$ is nonzero. Then if $S_1,S_2,\ldots,S_n$ are subsets of $\mathbb{F}$ with $|S_i|\ge t_i+1$, then there are $s_1\in S_1$,$s_2\in S_2$,\ldots,$s_n\in S_n$ so that $f(s_1,s_2,\ldots,s_n)\neq 0$.
\end{Thm}

In particular, a graph polynomial $P_{G}(\bm{x})$ is a homogeneous polynomial and $\deg(P_{G})$ is equal to $|E(G)|$. Therefore, if there exists a monomial $c\prod_{v \in V(G)} {x_{v}}^{t_{v}}$ in the expansion of $P_{G}$ so that $c \neq 0$ and $t_{v} < k$ for each $v \in V(G)$, then $G$ is $k$-choosable.
Jensen and Toft \cite{Jensen1995} defined the \emph{Alon-Tarsi number} of a graph as follows. 

\begin{Def}\textup{
The \emph{Alon-Tarsi number} of a graph $G$, denoted by $AT(G)$, is the minimum $k$ for which  there exists a monomial $c\prod_{v\in V(G)}x_v^{t_v}$ in the expansion of $P_{G}(\bm{x})$ such that $c\neq 0$ and $t_v<k$ for all $v\in V(G)$.}
\end{Def}

As explained above, $ch(G)\leq AT(G)$ for every graph $G$.
Moreover, it is known that the gap between $ch(G)$ and $AT(G)$ can be arbitrary large.
Nevertheless, it is also known that the upper bounds of $ch(G)$ and $AT(G)$ are the same for several graph classes.
For example, Thomassen \cite{Thomassen1994} proved that every planar graph is $5$-choosable. Later, Zhu proved the following.

\begin{Thm}[\cite{Zhu2018}]\label{Zhu2018}
Let $G$ be a plane graph. Then $AT(G)\leq 5$.
\end{Thm}

Moreover, it was shown in \cite{Cushing2010} that every planar graph is $1$-defective $4$-choosable. Recently, Grytczuk and Zhu have proved the following theorem.

\begin{Thm}[\cite{Grytczuk2019}]\label{Grytczuk2019}
Let $G$ be a plane graph. Then there exists a matching $M$ of $G$ such that $AT(G-M)\leq 4$.
\end{Thm} 
This result implies that every planar graph is $1$-defective $4$-choosable.
Furthermore, it was shown independently in \cite{Eaton1999} and \cite{Skrekovski1999} that every planar graph is $2$-defective $3$-choosable.
In this context, it seems natural to ask whether there exists a subgraph $H$ of $G$ such that $AT(G-E(H))\leq 3$ and $d_{H}(v)\leq 2$ for every $v\in V(G)$.
Since if it was true, this implies that every planar graph is $2$-defective $3$-choosable.
However,this is not true and it was shown in \cite{Kim2019} that there exists a planar graph $G$ such that for any subgraph of $H$ of $G$ with maximum degree at most $3$, $G-E(H)$ is not $3$-choosable.
On the other hand, the following was also proved in the same paper.
\begin{Thm}[\cite{Kim2019}]\label{Kim2019}
Let $G$ be a plane graph.
Then there exists a forest $F$ in $G$ such that $AT(G-E(F))\leq 3$. 
\end{Thm}

A graph $H$ is a {\it minor} of a connected graph $G$ if we obtain $H$ from $G$ by deleting or contracting some edges recursively.
A graph $G$ is {\it $H$-minor-free} if $H$ is not a minor of $G$.
If multiple edges appear by a contraction, we replace them with simple edge.

As another extension of Thomassen's result, it was shown in \cite{He2008} and \cite{Skrekovski1998} that every $K_5$-minor-free graph is $5$-choosable.
Moreover, it is also shown in \cite{Woodall2011} that every $K_5$-minor-free graph is $1$-defective $4$-choosable.
In this paper, we extend these results from choice number to Alon-Tarsi number. 

\begin{Thm}\label{main}
Let $G$ be a $K_5$-minor-free graph.
Then all of the following hold.
\begin{itemize}
\item[(i)] $AT(G)\leq 5$.
\item[(ii)] There exists a matching $M$ of $G$ such that $AT(G-M)\leq 4$.
\end{itemize}
\end{Thm}

\begin{Thm}\label{main-forest}
For every $K_5$-minor-free graph $G$, there exists a forest $F$ such that 
$G - E(F)$ is $2$-degenerate.
\end{Thm}

Thus we have the following corollary.

\begin{Cor}\label{main-Cor}
For every $K_5$-minor-free graph $G$, there exists a forest $F$ such that 
$AT(G - E(F)) \leq 3$.
\end{Cor}

This paper is organized as follows. In Section 2, we prepare some lemmas in order to show the main theorems. And in Section 3, we prove Theorem \ref{main} and Theorem \ref{main-forest}.
In Section 4, we have some remarks that Theorem \ref{main} and Corollary \ref{main-Cor} can be extended to singed graphs.


\section{Orientations and Alon-Tarsi number}

\subsection{An alternative definition of the Alon-Tarsi number.}
Indeed Alon-Tarsi number is already defined algebraically in Section 1, Alon and Tarsi \cite{Alon1992} found a combinatorial interpretation of the coefficient for each monomial in the graph polynomials in terms of orientations and Eulerian subgraphs. 
For an orientation $D$ of $G$, {\it $d_{D}^{+}(v)$} (resp. {\it $d_{D}^{-}(v)$}) denotes out-degree (resp. in-degree) of a vertex $v$ in $D$.
The maximum out-degree of $D$ is denoted by $\Delta^{+}(D)$.
A subgraph $H$ of $D$ is called {\it Eulerian} if $V(H)=V(G)$ and $d_{H}^{-}(v)=d_{H}^{+}(v)$ for every $v\in V(H)$ with respect to $D$.
Note that $H$ might not be connected.
Let {\it $EE(D)$} (resp. {\it $OE(D)$}) denote the set of all Eulerian subgraphs of $D$ with even (resp. odd) number of edges.
Especially, we say that an orientation $D$ is {\it acyclic} if $D$ does not contain any directed cycles.

\begin{Thm}[\cite{Alon1992}]\label{AT}
Let $G$ be a graph, let $P_{G}$ be the graph polynomial of $G$ and let $D$ be an orientation of $G$ with out-degree sequence $\bm{d}=(d_{v})_{v\in V(G)}$. Then the coefficient of $\prod_{v\in V(G)}{x_{v}^{d_{v}}}$ in the expansion of $P_{G}$ is equal to $\pm (|EE(D)|-|OE(D)|)$ 
\end{Thm}

We say that orientation $D$ of $G$ is an {\it AT-orientation} if $D$ satisfies $|EE(D)|-|OE(D)|\neq 0$.

\subsection{Orientations of planar graphs}
Now let us focus on planar graphs.
We say a plane graph $G$ is a {\it near triangulation} if each internal face in $G$ is triangular.
In the papers \cite{Grytczuk2019}, \cite{Kim2019} and \cite{Zhu2018}, the following are shown respectively.

\begin{Lem}\label{lem1}
Let $G$ be a plane graph with simple boundary cycle $C=v_{1}v_{2}...v_{m}$.
Then all of the following hold.
\begin{itemize}
\item[(i)]{\rm (\cite{Zhu2018})} $G$ has an AT-orientation $D$ such that $d^{+}_{D}(v_{1})=0$, $d^{+}_{D}(v_{2})=1$, $d^{+}_{D}(v_{i})\leq 2$ for each $i\in \{3,...,m\}$ and $d^{+}_{D}(u)\leq 4$ for each interior vertex $u$.
\item[(ii)]{\rm (\cite{Grytczuk2019})} There exists a matching $M$ and an AT-orientation $D$ of $G-M$ such that
$d^{+}_{D}(v_{1})=d^{+}_{D}(v_{2})=0$, $d^{+}_{D}(v_{i})\leq 2-d_{M}(v_i)$ for each $i\in \{3,...,m\}$ and $d^{+}_{D}(u)\leq 3$ for each interior vertex $u$.
\item[(iii)]{\rm (\cite{Kim2019})} There exists a forest $F$ in $G$ and an  acylic orientation $D$ of $G-E(F)$ such that
$d^{+}_{D}(v_{1})=d^{+}_{D}(v_{2})=0$, $d^{+}_{D}(v_{i})=1$ for each $i\in \{3,...,m\}$ and $d^{+}_{D}(u)\leq 2$ for each interior vertex $u$.
\end{itemize}
\end{Lem}

In order to show the main theorem, we need an orientation which has some stronger properties.

\begin{Lem}\label{lem4}
Let $G$ be a plane graph with a boundary cycle $v_1v_2v_3$.
Then all of the following hold.
\begin{itemize} 
\item[(i)]There exists a matching $M$ of $G$ and an AT-orientation $D$ of $G-M$ such that $M$ does not cover $v_3$, $d^{+}_{D}(v_1)=d^{+}_{D}(v_2)=0$, $d^{+}_{D}(v_3)=2$ and $d^{+}_{D}(y)\leq 3$ for $y\in V(G)-\{v_1,v_2,v_3\}$.
\item[(ii)]There exists a forest $F$ of $G$ and an acyclic orientation $D$ of $G-E(F)$ such that $v_1v_3\not\in E(F)$, $d^{+}_{D}(v_1)=d^{+}_{D}(v_2)=0$, $d^{+}_{D}(v_3)=1$ and $d^{+}_{D}(y)\leq 2$ for $y\in V(G)-\{v_1,v_2,v_3\}$.
\end{itemize}
\end{Lem}

\begin{figure}[htbp]
\begin{center}
\includegraphics[scale=0.5]{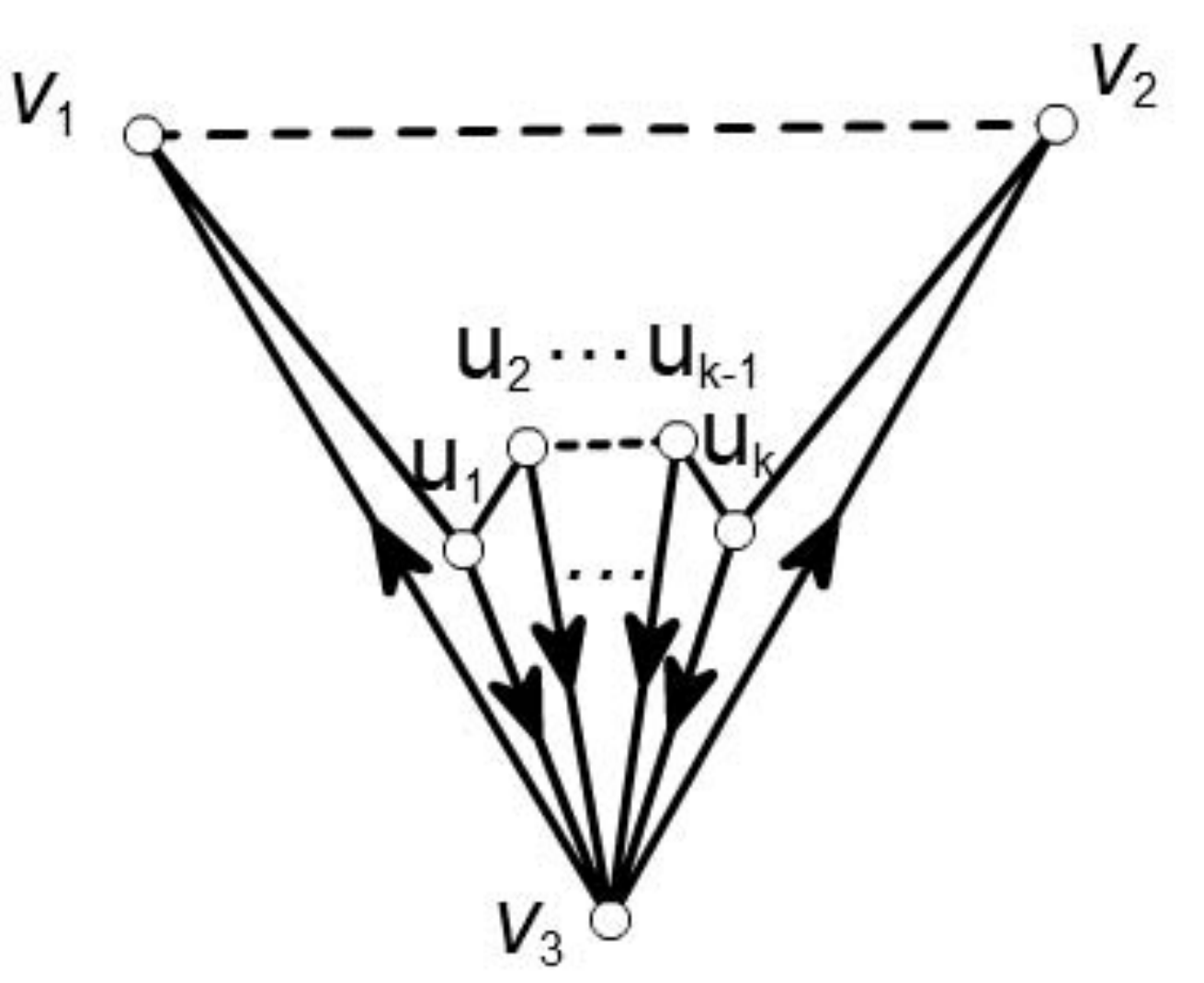}
\caption{The orientation $D$ of $G-M$}
\label{plane1}
\end{center}
\end{figure}

\noindent {\bf Proof.} Let $G'=G-v_{3}$ and let $N(v_3)=\{v_1,u_1,...,u_k,v_2\}$ be the neighborhood of $v_3$ as this rotation.
Since $G'$ is a plane graph, we have a matching $M$ and an AT-orientation $D'$ of $G'-M$ such that $d^{+}_{D'}(v_{1})=d^{+}_{D'}(v_{2})=0$, $d^{+}_{D'}(u_{i})\leq 2$ for $i\in \{1,2,...,k\}$ and $d^{+}_{D'}(u)\leq 3$ for each interior vertex $u$ by Lemma \ref{lem1} (See Figure \ref{plane1}).
Let $D$ be the orientation of $G-M$ obtained from $D'$ by adding the vertex $v_3$ and $k+2$ oriented edges $(u_i,v_3)$ for $i\in \{1,2,...,k\}$, $(v_3,v_1)$ and $(v_3,v_2)$.

It is easy to see that $D$ also satisfies the out-degree conditions and that $M$ does not cover $v_3$.
Moreover, since the vertices $v_1$ and $v_2$ have out-degree $0$, $D$ is also an AT-orientation. $\square$

\medskip
Let $G$ and $H$ be a graph which contain a clique of the same size. 
The {\it clique-sum} of $G$ and $H$ is a operation that forms a new graph obtained from their disjoint union by identifying a clique of $G$ and one of $H$ with the same size and possibly deleting some edges in the clique.
A {\it $k$-clique-sum} is a clique-sum in which both cliques have at most $k$-vertices.

\begin{Lem}\label{lem10}
Let $G$ be a graph which can be obtained by the $3$-clique-sum of $G_{1}$ and $G_{2}$ and let $T=\{x_1, x_2, x_3 \}$ be its clique.
Moreover, let $G'_{i}=G\cap G_{i}$.
Suppose that $G'_{1}$ has an AT-orientation $D'_{1}$ with $\Delta(D'_1)\leq k$ and that $G_{2}$ has an AT-orientation $D_{2}$ such that $d^{+}_{D_2}(x_1)=0$, $d^{+}_{D_2}(x_2)\leq 1$, $d^{+}_{D_2}(x_3)\leq 2$ and $\Delta(D_2)\leq k$ and $x_{i}$ is directed only to $x_{i'}$, where $x_{i},x_{i'}\in V(T)$.
Then $G$ has an AT-orientation $D$ such that $d^{+}_{D}(v)= d^{+}_{D'_{1}}(v)$ for each $v\in V(G_{1})$ and maximum out-degree of $D$ is at most $k$.
\end{Lem}

\noindent {\bf Proof.} Let $D'_{2} \subset D_{2}$ be the orientation of $G'_{2}$ and let $D=D'_{1}\cup \big(D'_{2}-E(T)\big)$.
Then it is easy to see that  $d^{+}_{D}(v)= d^{+}_{D'_{1}}(v)$ for each $v\in V(G'_{1})$, $d^{+}_{D}(v)= d^{+}_{D'_{2}}(v)$ for each $v\in V(G'_{2})-V(T)$ and hence maximum out-degree of $D$ is at most $k$.
For the orientation $D_{2}$, the vertices in $T$ has a direction only to other vertices of $T$ and no Eulerian subgraphs in $D_{2}$ contain the edge in $T$ by the out-degree conditions of $D_{2}$.
Thus $D'_{2}$ is also an AT-orientation of $G'_{2}$ and any spanning Eulerian sub-digraphs $H$ of $D$ has an edge-disjoint decomposition $H=H_{1}\cup H_{2}$ where $H_{1}$ and $H_{2}$ are Eulerian sub-digraphs in $D'_{1}$ and $D'_{2}$, respectively.
Therefore, we have the bijection $\tau$ so that 
\begin{itemize}
\item $\tau \big(EE(D)\big)=\big(EE(D'_{1}) \times EE(D'_{2})\big) \cup \big(OE(D'_{1})\times OE(D'_{2})\big)$ and
\item $\tau \big(OE(D)\big) = \big(OE(D'_{1})\times EE(D'_{2})\big) \cup \big(EE(D'_{1})\times OE(D'_{2}) \big)$.\\
\end{itemize}
Hence 
\begin{eqnarray*}\label{diffe}
&&|EE(D)|-|OE(D)|\\
&=&(|EE(D'_{1})|\times |EE(D'_{2})|+|OE(D'_{1})|\times |OE(D'_{2})|)\\
&&-(|EE(D'_{1})|\times |OE(D'_{2})|+|OE(D'_{1})|\times |EE(D'_{2})|)\\
&=&(|EE(D'_{1})|-|OE(D'_{1})|)\cdot(|EE(D'_{2})|-|OE(D'_{2})|)\\
&\neq&0.
\end{eqnarray*}
These imply that the orientation $D$ is an AT-orientation of $G$ with desired properties.  $\square$

\subsection{Characterizations of $K_{5}$-minor-free graphs.}

Now, let us focus on $K_{5}$-minor-free graphs.
In order to show the main theorem, we use the following results.

\begin{figure}[htbp]
\begin{center}
\includegraphics[scale=0.5]{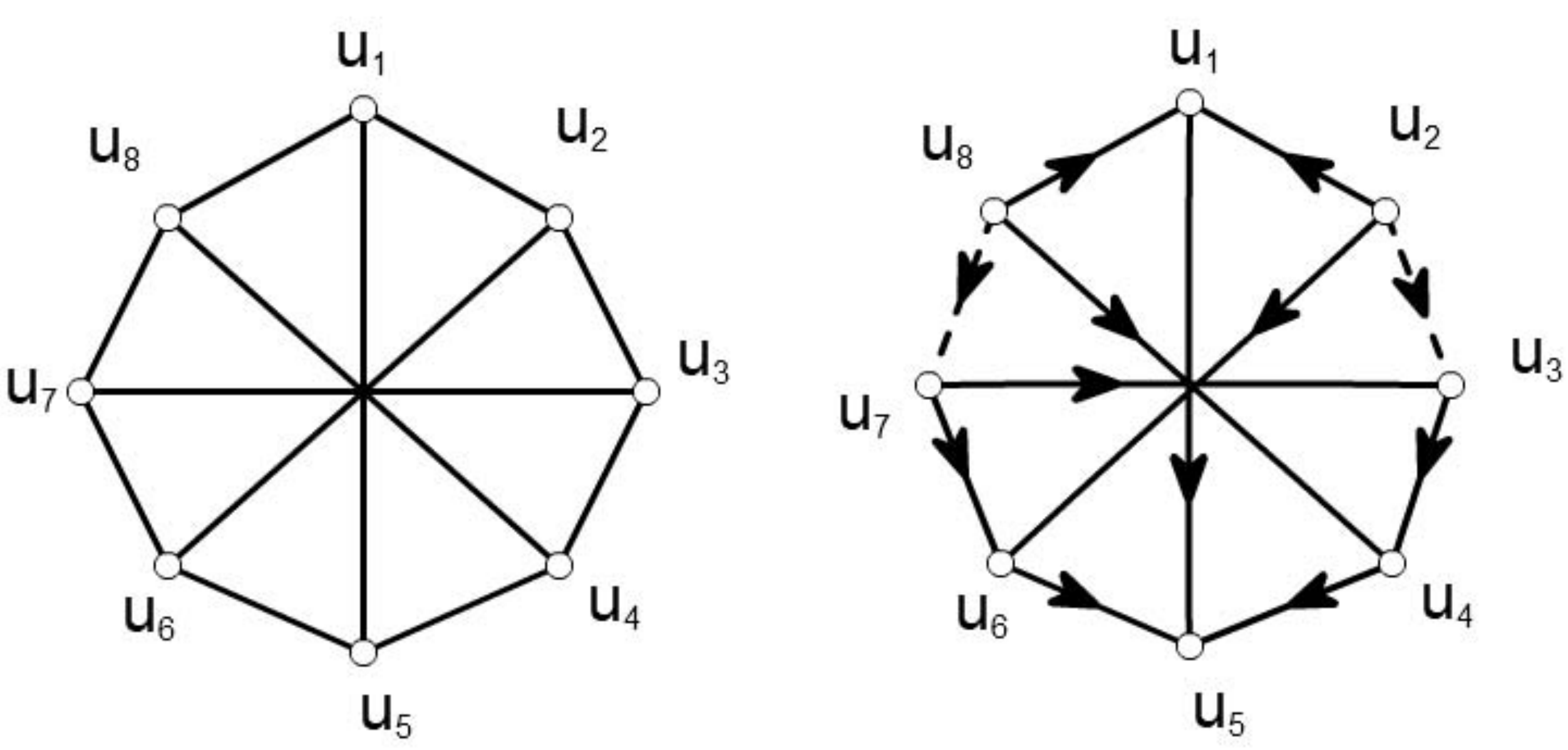}
\caption{The left is the Wagner graph $W$ and the right is an acyclic orientation with maximum out-degree $3$.
Doted edges denote elements of a matching or forest.
If we delete the two doted edges, the orientation has maximum out-degree $2$.}
\label{Wag}
\end{center}
\end{figure}

\begin{Lem}[\cite{Wagner1937}]\label{lem8}
A graph $G$ is $K_{5}$-minor-free if and only if $G$ can be formed from some $3$-clique-sums of graphs, each of which is either planar or the Wagner graph $W$ as shown in Figure \ref{Wag}.
\end{Lem}

\section{Proof of main Theorem.}

Theorem \ref{main} and Theorem \ref{main-forest} follow from the lemma below.

\begin{Lem}\label{lem2}
Let $G$ be a $K_5$-minor-free graph and let $H_{i}$ be a subgraph of $G$ which is isomorphic to $uv\in E(G)$ or $\{uv,vw,wu\}\subset E(G)$ for each $i\in\{1,2,3\}$.
Then all of the following hold.
\begin{itemize}
\item[(i)] There exists an AT-orientation $D$ such that $d^{+}_{D}(u)=0$, $d^{+}_{D}(v)=1$, ($d^{+}_{D}(w)=2$ if $H_1$ is isomorphic to $K_3$) and $d^{+}_{D}(y)\leq 4$ for $y\in V(G)-\{u,v,w\}$.
\item[(ii)]There exists a matching $M$ of $G$ and an AT-orientation $D$ of $G-M$ such that $d^{+}_{D}(u)=d^{+}_{D}(v)=0$, ($d^{+}_{D}(w)=2$ and $M$ does not cover $w$ if $H_2$ is isomorphic to $K_3$) and $d^{+}_{D}(y)\leq 3$ for $y\in V(G)-\{u,v,w\}$.
\item[(iii)]There exists a forest $F$ of $G$ and an acyclic orientation $D$ of $G-E(F)$ such that $uw\not\in E(F)$, $d^{+}_{D}(u)=d^{+}_{D}(v)=0$, ($d^{+}_{D}(w)=1$ if $H_3$ is isomorphic to $K_3$) and $d^{+}_{D}(y)\leq 2$ for $y\in V(G)-\{u,v,w\}$.
\end{itemize}
\end{Lem}

\noindent {\bf Proof.} Suppose that the Lemma is false and let $G_{i}$ be a counterexample for each $i\in \{1,2,3\}$ respectively with $|V(G_{i})|$ as small as possible.
By the minimality of $G_{i}$, $G_{i}$ is connected.
Moreover, it is easy to check that $G_{i}$ does not have a cut vertex.
Thus we may assume $G_{i}$ is $2$-connected.\\

First we suppose that $G_{i}$ is a plane graph.
If $H_{i}$ is isomorphic to $K_2$ or $K_{3}$ which bounds a face, without loss of generality, $H_{i}$ lies on the boundary of $G$.
In this case, a desired AT-orientation exists by Lemma \ref{lem1}.
Thus $H_{i}$ consists a separating $3$-cycle in $G_{i}$.
We let $G_{i,1}$ and $G_{i,2}$ be subgraphs of $G_{i}$ so that $G_{i,1}\cup G_{i,2}=G_{i}$ and $G_{i,1}\cap G_{i,2} = H_{i}$.
By Lemma \ref{lem1} and Lemma \ref{lem4}, for $j\in \{1,2\}$ we have an AT-orientation $D_{1,j}$ of $G_{1,j}$ which satisfies the conditions.
Similarly, we have that there exists a matching $M_{j}$ and an AT-orientation of $G_{2,j}-M_{j}$ and that there exists a forest $F_{j}$ and acyclic orientation of $G_{3,j}-E(F_{j})$, which satisfy the conditions.
It is easy to see that $M=M_{1}\cup M_{2}$ is also a matching of $G_{2}$ and $F=F_{1}\cup F_{2}$ is a forest of $G_{3}$.
Therefore, we get a desired AT-orientation respectively by Lemma \ref{lem10}.
Thus $G_{i}$ is not planar.

Next, suppose that $G_{i}$ is the Wagner graph $W$.
Since $W$ does not contain a triangle, $H_{i}$ must be isomorphic to $K_{2}$.
By the symmetry of $W$, we may assume that $H_{i}=u_1u_5$ or $u_5u_6$.
We let $M$ be a matching of $W$ and let $F$ be a forest of $W$ such that $M=E(F)=E(H)\cup \{u_2u_3, u_7u_8\}$.
In this case, the orientation in Figure \ref{Wag} is a desired AT-orientation respectively.

Thus we assume that $G_{i}$ is neither planar graph nor the graph $W$.
By Lemma \ref{lem8}, there exists $K_{5}$-minor-free graphs $G_{i,1}$ and $G_{i,2}$ such that $G_{i}$ can be obtained by a $3$-clique-sum of $G_{i,1}$ and $G_{i,2}$.
Let $T$ be its clique and let $G'_{i,j}=G_{i,j}\cap G_{i}$ for $j \in \{1,2\}$.
It is easy to see that $H_{i}\subset G'_{i,1}$ or $G'_{i,2}$.
Without loss of generality, we may assume that $H_{i}\subset G_{i,1}$.
By the minimality of $G_{i}$, we have the following.

\begin{itemize}
\item[(i)] There exists an AT-orientation $D_{1,1}$ of $G'_{1,1}$ which satisfies the assumption (i) of Lemma \ref{lem2}.
\item[(ii)] There exists a matching $M_1$ and an AT-orientation $D_{2,1}$ of $G'_{2,1}-M$ which satisfies the assumption (ii) of Lemma \ref{lem2}.
\item[(iii)] There exists a forest $F_1$ and an acyclic orientation $D_{3,1}$ of $G'_{3,1}-E(F)$ which satisfies the assumption (iii) of Lemma \ref{lem2}.
\end{itemize}
First, we consider the case when $i=1$.
By the minimality of $G_{1}$, we get an AT-orientation $D_{1,2}$ of $G_{1,2}$ with $d^{+}_{D_{1,2}}(x_1)=0$, $d^{+}_{D_{1,2}}(x_2)=1$, $d^{+}_{D_{1,2}}(x_3)=2$ and the maximum degree of $D_{1,2}$ is at most $4$.
By Lemma \ref{lem10}, we get a desired AT-orientation $D$ in $G_{1}$.   

Next, we consider the case when $i=2$.
By the minimality of $G_2$, we get a matching $M_2$ of $G_{2,2}$ and an AT-orientation $D_{2,2}$ of $G_{2,2}-M_2$ such that $d^{+}_{D_{2,2}}(x_1)=0$, $d^{+}_{D_{2,2}}(x_2)=0$, $d^{+}_{D_{2,2}}(x_3)=2$, maximum out-degree of $D_{2,2}$ is at most $3$ and $M_2$ does not cover $x_3$.
Let $M=M_1\cup \big(M_2- \{x_1x_2\}\big)$.
It is easy to see that $M$ is a matching of $G$.
By Lemma \ref{lem10}, we get a desired AT-orientation $D$ in $G_{2}-M$.

Finally, we consider the case when $i=3$.
By the minimality of $G_3$, we get a forest $F_2$ of $G_{3,2}$ and an acyclic orientation $D_{3,2}$ of $G_{3,2}-E(F_2)$ with $d^{+}_{D_{3,2}}(x_1)=d^{+}_{D_{3,2}}(x_2)=0$, $d^{+}_{D_{3,2}}(x_3)=1$ and maximum out-degree of $D_{3,2}$ is at most $2$.
Let $F=F_1\cup \left(F_2- E(T)\right)$.
Similarly, we can show that $F$ is a forest and $D_3$ is an acyclic orientation of $G_3 - E(F)$ with desired properties.
This is a contradiction and we completes the proof. $\square$

\bigskip
\noindent {{\bf [Proof of Theorem \ref{main} and Theorem \ref{main-forest}] }\\
Theorem \ref{main} follows immediately from (i) and (ii) in Lemma \ref{lem2}.
For Theorem \ref{main-forest}, each $K_{5}$-minor-free graph $G$ has a forest $F$ and and an acyclic orientation $D$ of $G-E(F)$ with maximum out-degree at most $2$ by Lemma \ref{lem2}.
Since $G-E(F)$ is finite and $D$ is acyclic, there exists a vertex $v$ with $d^{-}_{D}(v)=0$ and hence the vertex $v$ has degree at most $2$ in $G-E(F)$. $\square$ 


\section{Some remarks}
A {\it signed graph} is a pair $(G,\sigma)$, where $G$ is a graph and $\sigma$ is a {\it signature} of $G$ which assigns to each edge $e=uv$ of $G$ a sign $\sigma_{uv}\in \{1,-1\}$. Let
$$N_k =
\begin{cases}
\{0, \pm 1, \dots , \pm q\} & \text{if $k=2q+1$ is an odd integer,} \\
\{\pm 1, \dots , \pm q\} & \text{if $k=2q$ is an even integer.}
\end{cases}
$$
Note that $|N_k|=k$ for each integer $k$.
A {\it proper coloring} of $(G, \sigma)$ is a mapping $c:V(G)\to N_k$ such that $c(x)\neq \sigma_{xy}c(y)$ for each edge $xy$.
The {\it chromatic number} $\chi(G,\sigma)$ of $(G,\sigma)$ is minimum integer $t$ such that there exists a proper coloring $c:V(G)\to N_t$.
The choice number $ch(G,\sigma)$ of $(G,\sigma)$ is minimum integer $k$ such that for every $k$-list assignment $L$, there exists a proper coloring $c$ of $(G,\sigma)$ so that $c(v)\in L(v)$ for every $v\in V(G)$.

Let $(G,\sigma)$ be a signed graph and let `$<$' be an arbitrary fixed ordering of the vertices of $(G,\sigma)$.
The \emph{singed graph polynomial} of $(G,\sigma)$ is defined as 
 $$P_{G,\sigma}(\bm{x})=\prod_{u\sim v,u<v}(x_u-\sigma_{uv}x_v),$$
 where $u\sim v$ means that $u$ and $v$ are adjacent, and $\bm{x}=(x_v)_{v\in V(G)}$ is a vector of $|V(G)|$ variables indexed by the vertices of $G$.
It is easy to see that a mapping $c:V(G) \to \mathbb{Z}$ is a proper coloring of $(G,\sigma)$ if and only if $P_{G,\sigma}(\bm{c}) \neq 0$, where $\bm{c} = \big(c(v) \big)_{v \in V(G)}$.
The Alon-Tarsi number of $(G,\sigma)$ is defined similarly and we have $\chi(G,\sigma)\leq ch(G,\sigma)\leq AT(G,\sigma)$.

Let $(G,\sigma)$ be a singed graph and let $D$ be an orientation of $(G,\sigma)$.
Let $\sigma EE(D)$ (resp. $\sigma OE(D)$) denote the set of all spanning Eulerian sub-digraphs of $D$ with even (resp. odd) number of positive edges on $\sigma$.
It was shown in \cite{Wang2019} that Theorem \ref{AT} can be extended to signed one as follows.

\begin{Thm}[\cite{Wang2019}]
Let $(G,\sigma)$ be signed graph, let $P_{G,\sigma}$ be the signed graph polynomial of $(G,\sigma)$ and let $D$ be an orientation of $(G,\sigma)$ with out-degree sequence $\bm{d}=(d_{v})_{v\in V(G)}$. Then the coefficient of $\prod_{v\in V(G)}{x_{v}^{d_{v}}}$ in the expansion of $P_{G,\sigma}$ is equal to $\pm (|\sigma EE(D)|-|\sigma OE(D)|)$.
\end{Thm}

We say that an orientation $D$ of $(G,\sigma)$ is a $\sigma$AT-orientation if $D$ satisfies $|\sigma EE(D)|-|\sigma OE(D)|\neq 0$.

Let us focus on planar graphs.
In the papers \cite{Grytczuk2019} and \cite{Wang2019}, the following are shown respectively.

\begin{Lem}\label{lem11}
Let $(G,\sigma)$ be a signed near triangulation and let $C=v_{1}v_{2}...v_{m}$ be the boundary cycle of $G$.
Then all of the following hold.
\begin{itemize}
\item[(i)] {\rm (\cite{Wang2019})} $G$ has a $\sigma$AT-orientation $D$ such that $d^{+}_{D}(v_{1})=0$, $d^{+}_{D}(v_{2})=1$, $d^{+}_{D}(v_{i})\leq 2$ for each $i\in \{3,...,m\}$ and $d^{+}_{D}(u)\leq 4$ for each interior vertex $u$.
\item[(ii)] {\rm (\cite{Grytczuk2019})} There exists a matching $M$ and a $\sigma$AT-orientation $D$ of $G-M$ such that
$d^{+}_{D}(v_{1})=d^{+}_{D}(v_{2})=0$, $d^{+}_{D}(v_{i})\leq 2-d_{M}(v_i)$ for each $i\in \{3,...,m\}$ and $d^{+}_{D}(u)\leq 3$ for each interior vertex $u$.
\end{itemize}
\end{Lem}

Although Lemma \ref{lem11} only deals with near triangulations in the paper \cite{Wang2019}, it is not hard to extend the graph class from near triangulations to planar graphs.
Moreover, since all the arguments of Lemmas in Section 2 and Lemma \ref{lem2} work even if we replace AT-orientations into $\sigma$AT-orientations, we have the following results.

\begin{Thm}\label{main2}
Let $(G,\sigma)$ be a signed $K_5$-minor-free graph.
Then all of the following hold.
\begin{itemize}
\item[(i)] $AT(G,\sigma)\leq 5$.
\item[(ii)] There exists a matching $M$ of $G$ such that $AT(G-M,\sigma)\leq 4$.
\item[(iii)] There exsits a forest $F$ in $G$ such that $AT(G-E(F),\sigma)\leq 3$. 
\end{itemize}
\end{Thm}

\end{document}